# How IGN (France) computed the so-called 'centre of gravity' of physical Europe in 1989 and 2004

Jean-François HANGOUËT

*Univ. Gustave Eiffel, IGN-ENSG, Champs-sur-Marne, France*

jean-francois.hangouet@ensg.eu

c/o École nationale des sciences géographiques
6-8 Avenue Blaise Pascal
Cité Descartes
Champs-sur-Marne
77455 Marne-la-Vallée
CEDEX 2
France

Dr. Jean-François Hangouët graduated as a survey engineer from the École nationale des sciences géographiques in 1993. His professional and academic experience includes map generalization, photogrammetry, spatial data quality, teaching.

# How IGN (France) computed the so-called 'centre of gravity' of physical Europe in 1989 and 2004

Abstract: The method used by senior geodetic engineer Jean-Georges Affholder to determine what can be termed as the 'centre of gravity' of physical Europe in 1989 and 2004 relies on mathematical formulae which, in their only published version, happen to be flawed with typographical errors that do not reflect Mr. Affholder's actual mathematical exactness. This short epistemological paper summarizes the major steps of Affholder's method, provides a corrected version of the general formulae, and briefly recalls some particulars of the specific determination of the centre of gravity of physical Europe.

**Introduction**

Senior geodetic engineer Jean-Georges Affholder, then working at the French Institut géographique national (IGN, renamed Institut national de l'information géographique et forestière in 2012) devised, and perfected in the late 1980s, a method for computing the 'centre of gravity' of any geographical area. His method he applied to a wide range of areas (metropolitan France, European Union, Eurozone…), as documented by Affholder himself (1991, 2003). Although he retired in 2004, his method is still being used, see *e.g.* the results computed and presented by Hangouët (2016).

The 'centre of gravity' is of a conventional nature here, being not the direct measurement of some physical reality, but the output of mathematical operations applied to geographical data selected to represent the outer limits of national, regional or physical areas. This convention however is not without utility: relying on a proven, systematic process, it makes it possible to compare the positions of the centres of different areas, or to follow the trajectory with time of the centre of one area (*e.g.* of an administrative unit when limits change).

Jean-Georges Affholder also applied his method to determine the 'centre of gravity' of physical Europe, twice, in 1989 and 2004. Beside a few considerations bearing on the 1989 computation given in a chapter he wrote for a collective book in French (Affholder 1991), Jean-Georges Affholder hasn't written any scientific paper on the specific determination of the geographical centre of physical Europe. This very chapter also happens to be the sole communication where the mathematical formulae of his general method were stated. Regrettably, the formulae as published in that book are affected with misleading typographical errors.

**Affholder's method, in general**

The general principle of Affholder's method is as follows:

(1) The boundaries of the geographical area of interest are projected on a reference oblate spheroid. The result is an 'ellipsoidal polygon', the sides of which are made of geodetic arcs joining vertices of known geographical (longitude, latitude) coordinates.

(2) This polygon is curved in space, looking like the fragment of an eggshell made of an infinitely thin yet homogeneous matter. Its centre of gravity, as defined in statics, is not to be found in the polygon itself, but floating somewhere in space (in that part of space that is on the hollow side of the cup). The three-dimensional Cartesian coordinates of the centre of gravity are computed by means of infinitesimal calculus (computation of the centres of gravity and of the areas of a great number of exceedingly narrow strips covering the polygon extent, and combination of these elementary centres and areas to construct the centre of gravity of the polygon proper).
(3) This point in space is projected orthogonally back onto the spheroid, providing the geographical coordinates of the looked-for centre of the polygon.

**Affholder's application to physical Europe**

Affholder applied his general method to different zones covering France and/or other European countries. He also applied it to physical Europe, twice, in 1989 and 2004.

For this application, the selection of the constituent parts of the so-called 'physical Europe' followed the extensive definition provided by Lapie (1838), with adjustments for some Atlantic, Mediterranean and Barents islands to account for modern evidence on continental shelves.

The boundaries of physical Europe were specifically constituted:

- in 1989, by digitizing IGN's 1 : 33 000 000 scale map of the World,
- in 2004, more precisely, from DCW/VMap 0 datasets, completed by the natural limits of physical Europe (crest lines, rivers…) that run through the Russian and Turkish territories, digitized from available 1 : 1 000 000 scale maps.

On both occasions, the Hayford ellipsoid was chosen as reference spheroid.

**Affholder's formulae**

In the general method described above, step 2 certainly makes the most specific problem (steps 1 and 3 involving but basic geodetic equations). Here are the mathematical formulae established by Jean-Georges Affholder to solve it, not as they were published, but in their corrected version.

In these formulae:

- a and b are the two fundamental parameters of the oblate spheroid
- $\lambda_0$ is an arbitrary reference longitude, the value of which will be constant
- $P_i(\lambda_i, \varphi_i)$ and $P_{i+1}(\lambda_{i+1}, \varphi_{+1i})$ are two successive vertices of the boundary of the ellipsoidal polygon, and $M(\lambda_M, \varphi_M)$ is their 'midpoint', defined as follows: $\lambda_M = (\lambda_i + \lambda_{i+1})/2$ and $\varphi_M = (\varphi_i + \varphi_{i+1})/2$
- $e^2 = \frac{a^2 - b^2}{a^2}$
- $N_M = \frac{a}{(1 - e^2 \cdot \sin^2 \varphi_M)^{1/2}}$
- $\rho_M = \frac{a(1 - e^2)}{(1 - e^2 \cdot \sin^2 \varphi_M)^{3/2}}$

Infinitesimal strips of spheroid are introduced (Fig. 1).

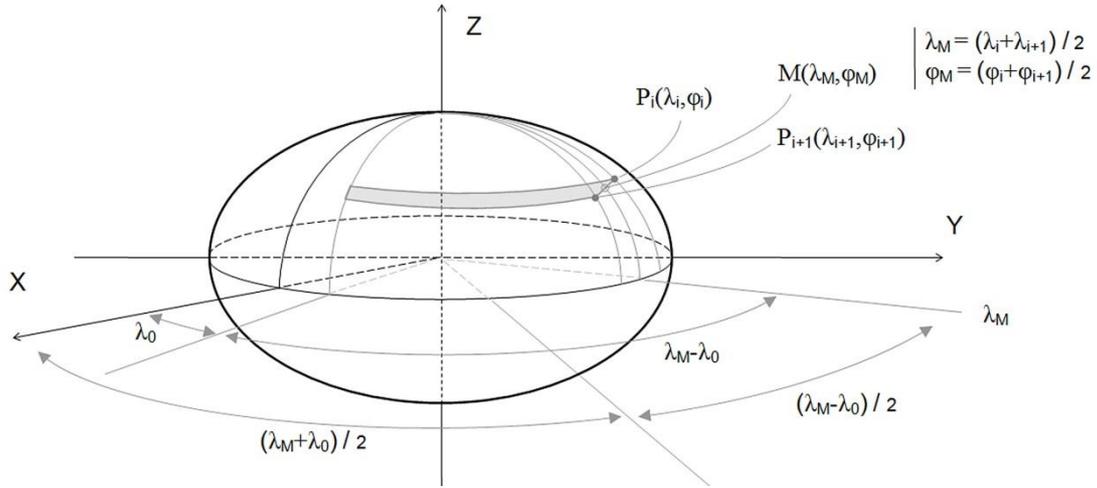

Figure 1. An infinitesimal strip on the oblate spheroid.

In practice, given the resolution of the input data, these strips are very small indeed: if extending, longitudinally, between $\lambda_0$ and $\lambda_M$, their latitudinal extent is $(\varphi_{i+1} - \varphi_i)$. They look very much like circular arcs. Their length and width can be approximated as, respectively:

$$\text{length}_i = N_M \cdot \cos \varphi_M \cdot (\lambda_M - \lambda_0)$$
$$\text{width}_i = \rho_M \cdot (\varphi_{i+1} - \varphi_i)$$

The centre of gravity $G_i = (X_i\ Y_i\ Z_i)$ of any such elementary strip can be considered to lie on the plane that passes through M and is parallel to the Z=0 plane. The elementary strip showing a circular shape, the distance $d$ between $G_i$ and the Z axis can be approximated as the distance of the centre of gravity of a circular arc to the centre of the circle. With R the radius of the circle, and $2\cdot\alpha$ the angle at the centre, this distance is: $d = R \cdot \frac{\sin \alpha}{\alpha}$, as demonstrated by Crawford (1898), among others. Here, $\alpha = \frac{\lambda_M - \lambda_0}{2}$ and $R = N_M \cdot \cos \varphi_M$.

With these elements, and some trigonometric rearrangements, the signed area $S_i$ and the centre of gravity $G_i = (X_i\ Y_i\ Z_i)$ of any such elementary strip are, with a first-order approximation:

$$S_i = N_M \cdot \cos \varphi_M \cdot (\lambda_M - \lambda_0) \cdot \rho_M \cdot (\varphi_{i+1} - \varphi_i)$$
$$X_i = N_M \cdot \cos \varphi_M \cdot \frac{\sin \lambda_M - \sin \lambda_0}{\lambda_M - \lambda_0}$$
$$Y_i = N_M \cdot \cos \varphi_M \cdot \frac{\cos \lambda_0 - \cos \lambda_M}{\lambda_M - \lambda_0}$$
$$Z_i = N_M \cdot (1 - e^2) \cdot \sin \varphi_M$$

The area S of the ellipsoidal polygon and the three-dimensional Cartesian coordinates of its centre of gravity can now be computed as, respectively, the total sum of the elementary areas and the weighted mean of the elementary coordinates:

$$S = \left|\sum_{i=1}^{n} S_i\right| \qquad X_G = \frac{\sum_{i=1}^{n} S_i \cdot X_i}{\sum_{i=1}^{n} S_i} \qquad Y_G = \frac{\sum_{i=1}^{n} S_i \cdot Y_i}{\sum_{i=1}^{n} S_i} \qquad Z_G = \frac{\sum_{i=1}^{n} S_i \cdot Z_i}{\sum_{i=1}^{n} S_i}$$

**Conclusion**

The two centres Jean-Georges Affholder found, differing from some 5 km, were consistent:

- 1989 determination:
  latitude:   54°54' N,         longitude:  25°19' E
- 2004 determination:
  latitude:   54°50'45" N,    longitude:  25°18'23" E

The latter, and more precise location, lies some 26 km north of Vilnius, in Lithuania, near the locality called Purnuškės. A monument was erected there in 2004, celebrating both the centre of gravity of physical Europe (in the geographical and mathematical sense reviewed here) and Lithuania's membership to the European Community, now European Union (in the political and economical sense exposed in the Maastricht Treaty).